\documentclass{amsart}
\usepackage{showkeys}
\usepackage{amsmath} 
\usepackage{amssymb}
\usepackage{mathrsfs}

\newtheorem{theorem}{Theorem}[section]

\newtheorem{lemma}[theorem]{Lemma}

\theoremstyle{definition}
\newtheorem{definition}[theorem]{Definition}

\theoremstyle{remark}
\newtheorem{remark}[theorem]{Remark}

\newtheorem{notation}[theorem]{Notation}
\newtheorem{context}[theorem]{Context}

\def\smallbox#1{\leavevmode\thinspace\hbox{\vrule\vtop{\vbox
   {\hrule\kern1pt\hbox{\vphantom{\tt/}\thinspace{\tt#1}\thinspace}}
   \kern1pt\hrule}\vrule}\thinspace}

\newcommand{\st}{{\rm such that}}

\newcommand{\ub}{{\rm ub}}

\newcommand{\Range}{{\rm Range}}

\newcommand{\Rang}{{\rm Rang}}

\newcommand{\then}{{\underline{then}}}
\newcommand{\when}{{\underline{when}}}
\newcommand{\Then}{{\underline{Then}}}

\newcommand{\mn}{{\medskip\noindent}}
\newcommand{\sn}{{\smallskip\noindent}}

\newcommand{\cD}{{\mathscr D}}

\newcommand{\cH}{{\mathscr H}}

\newcommand{\cE}{{\mathscr E}}

\newcommand{\cM}{{\mathscr M}}

\newcommand{\cP}{{\mathscr P}}

\newcommand{\bbZ}{{\mathbb Z}}

\newcommand{\cS}{{\mathscr S}}

\newcommand{\cU}{{\mathscr U}}

\newcommand{\cf}{{\rm cf}}

\newcount\skewfactor
\def\mathunderaccent#1#2 {\let\theaccent#1\skewfactor#2
\mathpalette\putaccentunder}
\def\putaccentunder#1#2{\oalign{$#1#2$\crcr\hidewidth
\vbox to.2ex{\hbox{$#1\skew\skewfactor\theaccent{}$}\vss}\hidewidth}}
\def\name{\mathunderaccent\tilde-3 }
\def\Name{\mathunderaccent\widetilde-3 }

\newenvironment{PROOF}[2][\proofname.]
   {\begin{proof}[#1]\renewcommand{\qedsymbol}{\textsquare$_{\rm #2}$}}
   {\end{proof}}

\begin{document}

\newcount\skewfactor
\def\mathunderaccent#1#2 {\let\theaccent#1\skewfactor#2
\mathpalette\putaccentunder}
\def\putaccentunder#1#2{\oalign{$#1#2$\crcr\hidewidth
\vbox to.2ex{\hbox{$#1\skew\skewfactor\theaccent{}$}\vss}\hidewidth}}
\def\name{\mathunderaccent\tilde-3 }
\def\Name{\mathunderaccent\widetilde-3 }

\newcommand{\E}{{\varepsilon}}
\newcommand{\eps}{{\epsilon}}
\newcommand{\f}{{\varphi}}
\newcommand{\tet}{{\theta}}
\newcommand{\pS}{\subset}
\newcommand{\ps}{\subseteq}
\newcommand{\Ps}{\nsubseteq}
\newcommand{\pa}{\forall}
\newcommand{\py}{\exists}
\newcommand{\pg}{\models}
\newcommand{\pc}{\wedge}
\newcommand{\pd}{\vee}
\newcommand{\po}{\emptyset}
\newcommand{\lk}{\langle}
\newcommand{\rk}{\rangle}
\newcommand{\pL}{\parallel}
\newcommand{\up}{\Gamma}
\newcommand{\ha}{{\aleph}}
\newcommand{\wto}{\Rightarrow}
\newcommand{\uph}{\upharpoonright}

\newcommand{\gK}{\mathfrak{K}}
\newcommand{\gq}{\mathfrak{q}}
\newcommand{\gJ}{\mathfrak{J}}
\newcommand{\gZ}{\mathfrak{Z}}

\newcommand{\high}{{\rm high}}
\newcommand{\cut}{{\rm cut}}
\newcommand{\sau}{{\rm sau}}
\newcommand{\siu}{{\rm siu}}
\newcommand{\ext}{{\rm Ext}}
\newcommand{\Depth}{{\rm Depth}}

\newcommand{\multi}{{multiplicatively}}
\newcommand{\seq}{{sequence}}
\newcommand{\cont}{{continuous}}
\newcommand{\incr}{{increasing}}
\newcommand{\Def}{{Definition}}
\newcommand{\homo}{{homomorphism}}
\newcommand{\stat}{{stationary}}
\newcommand{\Wlog}{{without loss of generality}}


\setcounter{section}{-1}


\title{Compactness in singular cardinals revisited}
\author{Saharon Shelah}
\address{Institute of Mathematics
 The Hebrew University of Jerusalem
 Jerusalem 91904, Israel
 and  Department of Mathematics
 Rutgers University
 New Brunswick, NJ 08854, USA}
\email{shelah@math.huji.ac.il}
\urladdr{http://www.math.rutgers.edu/\char`\~shelah}
\thanks{Research supported by the United States-Israel Binational
Science Foundation.  In a mimeographed from this was included
in May 1977, and a lecture on it was given in Berlin 1977.  
References like \cite[Th0.2=Ly5]{Sh:950} means the label of Th.0.2 is y5.
The reader should note that the version in my website is usually more updated
than the one in the mathematical archive.
Publication 266}




 
\date{January 15, 2019}

\begin{abstract}
This is the second combinatorial proof of the compactness theorem for
singular from 1977.  In fact it gives a somewhat stronger theorem.
\end{abstract}

\maketitle
\numberwithin{equation}{section}
\setcounter{section}{-1}
\newpage

\section {Introduction}

For a long time I have been interested in compactness in singular cardinals;
i.e., whether if something occurs for ``many" subsets of a singular $\lambda$
 of cardinality $<\lambda$, it occurs for $\lambda$. For the positive side
in the seventies we have

\begin{theorem}
\label{0.1}
Let $\lambda$ be a singular cardinal, $\chi^*<\lambda$.
Let $\cU$ be a set, ${\bf F}$ a family
of pairs $(A,B)$ of subsets of ${\cU}$, instead of $(A,B) \in {\bf F}$
we may write $A/B \in {\bf F}$ (formal quotient) or $A/B$
is ${\bf F}$-free. Assume further that ${\bf F}$ is a nice freeness notion
meaning it satisfies axioms II, III, IV,VI, VII from \ref{0.2} below. 
Let $A^*, B^* \subseteq {\cU}$ with $|B^*|=\lambda$. 

\Then \, $B^* / A^*\in {\bf F}$ is free in a weak sense, that is: there is
an increasing continuous sequence $\langle A_\alpha:\alpha <
\delta\rangle$ of subsets of $B_*$ of cardinality $< \lambda$ such
that $A_0 = \emptyset,\bigcup\limits_{\alpha < \delta} A_\alpha = A_*$
and $A_{i+1}/A_i \cup A$ is $\bold F$-free for $i < \lambda$ 
\when \, (see Definition \ref{0.3} below):
\mn
\begin{enumerate}
\item[$(*)_0$] for the ${\cD}_{\chi^*}(B^*)$-majority of
$B \in [B^*]^{<\lambda}$ we have $B/ A^* \in {\bf F}$

\underline{or just}
\sn
\item[$(*)_1$]  the set $\{\mu < \lambda:\{B \in [B^*]^\mu:B/ A^* \in
{\bf F}\} \in \cE^{\mu^+}_\mu (B^*)\}$ contains a club of $\lambda$,

\underline{or at least}
\sn
\item[$(*)_2$]  for some set $C$ of cardinals $<\lambda$, 
unbounded in $\lambda$  and closed (meaningful only if 
$\cf(\lambda)>\aleph_0)$, for every $\mu \in C$,
for an ${\cE}^{\mu^+}_\mu (B^+)$-positive set of $B \in [B^*]^\mu$ 
we have $B/ A^*\in{\bf F}$.
\end{enumerate}
\end{theorem}

\noindent
Where
\begin{definition}
\label{0.2}
For a set ${\cU}$ and ${\bf F} \subseteq \{(A,B):A,B \subseteq {\cU}\}$
we say, ${\bf F}$ is a $\chi$-nice freeness notion \underline{if} ${\bf F}$
satisfies:
\mn
\begin{enumerate}
\item[Ax.II]  $B/ A \in {\bf F} \Leftrightarrow A \cup  B/A \in {\bf F}$
\sn
\item[Ax.III]  if $A \subseteq B \subseteq C, B/A \in
{\bf F}$ and $C/B \in {\bf F}$ then $C/ A\in {\bf F}$,
\sn
\item[Ax.IV]  if $\langle A_i:i \le \theta\rangle$ is increasing
continuous, $\theta = \cf(\theta),A_{i+1}/A_i \in {\bf F}$ then
$A_\theta/A_0 \in{\bf F}$,
\sn
\item[Ax.VI]  if $A/B \in {\bf F}$ then for the
${\cD}_{\chi}$-majority of $A'\subseteq A$ we have, $A'/B \in {\bf F}$
(see below),
\sn
\item[Ax.VII]  if $A/B \in{\bf F}$ then for the
${\cD}_{\chi}$-majority of $A'\subseteq A$ we have, $A/B \cup A' \in
{\bf F}$.
\end{enumerate}
\end{definition}

\begin{definition}
\label{0.3}
1) Let ${\cD}$ be a function giving for any set $B^*$ a filter
${\cD}(B^*)$ on ${\cP}(B^*)$ (or on $[B^*]^\mu$).  \Then \, to say ``for the
${\cD}$-majority of $B \subseteq B^*$ (or $B \in [B^*]^\mu)$  we have
$\varphi(B)$" means $\{B \subseteq B^*:\varphi(B)\} \in \cD(B^*)$
(or $\{B \in[B^*]^\mu:\neg \varphi(B)\} = \emptyset \mod \cD(B^*))$.

\noindent
2) Let ${\cD}_\mu(B^*)$ be the family of $Y \subseteq {\cP}(B^*)$ 
such that for some algebra $M$ with universe $B^*$ and $\le \mu$ functions,
\[
Y \supseteq S_M = \{B \subseteq B^*:B \ne \emptyset 
\text{ is closed under the functions of } M\}.
\]

\noindent
2A) Let $\cD_{=\mu}(B^*)$ be defined similarly considering only $B$'s of
cardinality $\le \mu$.

\noindent
3) ${\cE}^{\mu}_\kappa (B^*)$ where $\mu \le \kappa^+$ is the collection of all
$Y \subseteq [B^*]^\kappa$ such that: for some $\chi,x$ satisfying $\{B^*,x\}
\in {\cH}(\chi)$, if $\bar{M} = \langle M_i:i<\mu\rangle$ is an
increasing continuous sequence of elementary submodels of $({\cH}(\chi),\in)$
such that $x \in M_0,\kappa +1 \subseteq M_0,\|M_i\| = \kappa$ 
and $i < \mu \Rightarrow \bar{M} \restriction (i+1)
\in M_{i+1}$, \then \, 
\mn
\begin{enumerate}
\item[$(a)$]  if $\mu \le \kappa$ then $\bigcup\limits_{i<\mu} M_i 
\cap B^*\in Y$
\sn
\item[$(b)$]  if $\mu=\kappa^+$ then for some club $C$ of $\mu^+$ we
 have $i \in C \Rightarrow M_i \cap B^*\in Y$.
\end{enumerate}
\end{definition}

\noindent
On ${\cD}_\mu$ see Kueker $[Ku]$, and on ${\cE}^{\mu^+}_\mu$
see \cite{Sh:52} repeated in \S2 below, note that in \cite{Sh:52} 
the axioms are phrased with elementary submodels rather then saying 
``majority".  The theorem was proved in \cite{Sh:52} but with two 
extra axioms, however it included the full case for varieties 
(i.e., including the non-Schreier ones). Later, the author eliminated
those two extra axioms: Ax. V and Ax. I. Now Ax. V was used in one point only
in \cite[\S1]{Sh:52}, and I eliminated it early (as presented in
\cite{BD}).  Axiom I is more interesting: it say that if $A'\subseteq A$ and
$A/B$ free then $A'/B$ is ${\bf F}$-free''; this is
like ``every subgroup of a free group if free; (this was shown not to
be necessary for varieties already in \cite{Sh:52}).
 
In 77 Fleissner has asked for a simpler ``combinatorial" proof and we find
such proof circulateding it in mimeographed notes \cite{Sh:E18}.  In 
May 77, and lecture on it in Berlin (summer 77 giving the full details
only for the case close to Abelian groups). This proof eliminates
the two extra axioms (as its assumptions holds by 
\cite[Lemma 3.4,p.349]{Sh:52}, see \S2 below).

Continuing this Hodges do \cite{Ho81} which contain a compactness result
and new important applications. I have thought he just represent the theorem
but looking at it lately it seems to me this is not exactly so; the
main point in the proof appears but the frame is different so it is
relative.  This exemplifies the old maxim ``if you want things 
done in the way you want it, you have to do them yourself". 

Anyhow below in \S1,\S2 we repeat the mimeographed notes.  Note that
\S2 repeats \cite[3.4]{Sh:52} needed for deducing \ref{0.1}.
Restricted to the needed case; note \ref{3.4} give hypothesis I (the non
${\cE}^{\lambda^+_i}_{\lambda^\gamma_i}$-non freeness is $(*)_2$ of
\ref{0.1} where hypothesis II is a weak form of Ax VII.

We thank Wilfred Hodges for help with some corrections and
encouragement and Paul Eklof for preserving and giving me a copy 
of the mimeographed notes after many years.
\newpage

\section {A compactness theorem for singular}

Here we somewhat improve and simplify the proof of \cite{Sh:52}
(and \cite{BD}).  It may be considered an answer to question B2 of
Fleissner \cite{Fli77}.

\begin{theorem}
\label{1.1}
Assume
\mn
\begin{enumerate}
\item[$(a)$]  $\lambda$ is a singular cardinal, $\lambda_i(i<\kappa)$ an
increasing and continuous sequence of cardinals (we let 
$\lambda(i)=\lambda_i)$ and
\[
\lambda_0=0, \kappa= \cf(\lambda), \kappa \le \lambda_1,
\lambda = \sum\limits_{i<\kappa} \lambda_i.
\]
\sn
\item[$(b)$] Let $S_i=\{A \subseteq \lambda: |A|=\lambda_i\}$ and 
$S'_i= S_i \cup \{\emptyset\}$
\sn
\item[$(c)$]  ${\bf F}$ is a family\footnote{Note that none of the
  axioms of \ref{0.2} is assumed.}
of pairs $(A,B), \lambda \supseteq A \supseteq B$; we may 
write ``$A/B$ belong to ${\bf F}$"
\sn
\item[$(d)$]  \underline{hypothesis I}: for each 
$i,i<\kappa, i$ a successor, there is a function $g_i$, two-place,
from $S'_i$ to $S'_i,$ such that: if $A_1 \subseteq A_2$ are from
$S'_i, A_1 \in \{\emptyset\} \cup \Range(g_i)$, \then \,
$A_2 \ps g_i (A_1, A_2)$ and $[g_i(A_1,A_2)/A_1] \in \bold F$
\sn
\item[$(e)$]  \underline{hypothesis II}:  if $i<\kappa,\, A,B \in 
S'_{i+1}$, $A\subseteq B$ and $B/A \in \bold F$ and $B \in
\Rang(g_{i+1})$, \then \, player II has a winning strategy in the 
following game $Gm_i [A,B]$. In the $n$-th move
$(n<\omega)$ player I choose $A_n \in S_i,$ such that $B_{n-1} \subseteq
A_n$, and then player II choose $B_n$, such that $A_n \subseteq 
B_n \in S_i$ \, (where we stipulate $B_{-1}= \emptyset)$.
Player II wins in the play if $(B \cup
\bigcup\limits_{n<\omega} B_n,A \cup \bigcup\limits_{n<\omega} B_n) 
\in {\bf F}$ (for $i=0$ this is an empty demand as
$S'_i=\{\emptyset\}$).
\end{enumerate}
\mn
\Then \, we can find an \incr\ and \cont\ chain
$A_\alpha (\alpha < \omega\kappa),$
\st\ $A_0=\emptyset, \lambda= \bigcup\limits_{\alpha} A_\alpha $ and
$A_{\alpha+1}/ A_\alpha \in {\bf F}$ for each $\alpha$.
\end{theorem}

\begin{PROOF}{\ref{1.1}}
Let in Hypothesis II the winning strategy of player II in the game
$Gm_i$ be given by the functions $h^n_i (A_0,\ldots,A_n; A,B)$. 
We define by induction on $i < \omega$ sets $A^n_i, B^n_i$ 
(for $i<\kappa)$ such that:
\mn
\begin{enumerate}
\item   $A^n_i (i<\kappa)$ is increasing and continuous in $i$ and
$A^n_i, B^n_i \in S_i$
\sn
\item   $A^n_i \subseteq B^n_i \subseteq A^{n+1}_i$
\sn
\item   $(B^n_i / B^{n-1}_i) \in {\bf F}$ where we stipulate
$B^{-1}_i = \emptyset$ and $B^n_i \in \Rang(g_i)$ for $i$ successor
\sn
\item   for $i < \kappa, 0 \le m<n, i<\omega$ we have
\[
h^{n-m}_i (A^{m+1}_i, A^{m+2}_i,\ldots,A^n_i; B^{m-1}_{i+1}, B^m_{i+1})
\subseteq A^{n+1}_i.
\]
\end{enumerate}
\mn
\underline{For $n=0$}

Let $A^0_i = \lambda_i, B^0_i = g_i (\emptyset,\lambda_i)$; clearly 
condition (1) holds, (2) and (4) say nothing and condition (3) 
holds by Hypothesis I.
\smallskip

\noindent
\underline{For $n+1$} assuming that for $n$ we have defined.

Let

\[
C^n_i = \bigcup\limits_{m< n} h^{n-m-1}_i (A^{m+1}_i,
A^{m+2}_i, \ldots, A^n_i, B^{m-1}_{i+1}, B^m_{i+1}) \cup B^n_i
\]

\mn
clearly $|C^n_i| = \lambda_i,$ hence we can let
$C_i^n= \{ {}^{\alpha}b^n_i :  \alpha <\lambda_i\}$.

Now we define $A^{n+1}_i = \{ {}^{\alpha}b^n_j : j<\kappa, \alpha
<\min \{\lambda_i,\lambda_j\}\}.$

Clearly condition (4) and the relevant parts of conditions (1) and (2) hold.
We have to choose $B^{n+1}_i$ such that

\[
A^{n+1}_i \subseteq B^{n+1}_i \text{ and } |B^{n+1}_i| =
\lambda_i, \text{ and } i \text{ successor } \Rightarrow
B^{n+1}_i/ B^n_i \in {\bf F}.
\]

\mn
So we let $B^{n+1}_i= g_i(B^n_i,A_i^{n+1})$ except that $B^{n+1}_0 =
\emptyset$.  By Hypothesis I this is O.K.

Now we can prove the conclusion of the theorem. We let
$D_{\omega i+k} = (B^{k-1}_{i+1} \cap \bigcup\limits_{m<\omega} A^m_i) \cup
\bigcup\limits_{j<i,m<\omega} A^m_j$ for
$i<\kappa$.   Clearly $D_0= \emptyset$, (in fact $A^n_i,B^n_i$ are $\emptyset$
for $i=0$); $\lambda= \bigcup\limits_{i<\omega \kappa} D_i$ as 
$\lambda_i = A^0_i \subseteq D_{\omega(i+1)} \subseteq \lambda$.  The sequence
is increasing and continuous. [that is e.g., if $\delta = \omega i + \omega$ so
$\delta=\omega (i+1)+0$ then $D_\delta \subseteq \bigcup\limits_{\alpha<\delta}
D_\alpha$ as $B^{-1}_{i+1} = \emptyset$, so $D_\delta=
\bigcup\limits_{j<i+1,m<\omega} A^m_j = 
(\bigcup\limits_{j<i,m<\omega} A^m_j) \cup 
\bigcup\limits_m A^m_i$ but $A^m_i \subseteq A^m_{i+1} \subseteq
B^m_{i+1}$ so $D_\delta \subseteq \bigcup\limits_{k}
\big[\bigcup\limits_{\binom{j<i}{m<\omega}} A^m_j \cup
(B^{k-1}_{i+1}\cap \bigcup\limits_{m<\omega} A^m_i) =
\bigcup\limits_k \, D_{\omega i+k} \subseteq \bigcup\limits_{\alpha<\delta}
D_\alpha \big]$.

Now $D_{\omega i+k+1} / D_{\omega i+k} \in {\bf F}$ as 
$B^k_{i+1} / B^{k-1}_{i+1} \in {\bf F}$ by condition (3),
 and then use condition (4)) and the choice of the $h^i_n -s'$
[that is, player II wins the play $\langle A^{k+\ell}_i, h^{k+\ell}_i
(A^{k+1}_i, A^{k+2}_i \ldots, A^{k+\ell}_i, B^{k-1}_{i+1}, B^k_{i+1}):
\ell<\omega\rangle$ of the game $Gm_i[B^{k-1}_{i+1}, B^k_{i+1}]$].
\end{PROOF}

\begin{remark}
\label{1.2}
1)  In the context of \cite{Sh:2}, \cite{BD} Hypothesis I holds
quite straightforwardly whereas Hypothesis II is proved separately, see
\cite[Lemma 3.4 p. 344]{Sh:52}.

\noindent
2)   Usually the choice of the $\lambda_i$'s is not important, and
then Hypothesis I, Hypothesis II should speak on
$\mu < \lambda,\mu < \mu'< \lambda$.

\noindent
3)  In the construction proving the Theorem we can continue $\chi <
\lambda_1$ steps instead of $\omega$ steps. We succeed if: in Hypothesis II
the game has length $\chi$ and we add to hypothesis II: if
$A_i / A_0 \in {\bf F}$ for $i< \chi,A_i$ increasing continuous 
then $\bigcup\limits_{i<\chi} A_i / A_0 \in {\bf F}$.

An example is: $G$ is a group with universe $\lambda$ and ${\bf F}=\{(A,B):
{\rm Ext} (A/B, C^*)= \emptyset\}$ where $A \subseteq B$ are subgroup of $G$,
$\cf(\lambda) < \chi < \lambda,\chi$ measurable $(C^*$ a fixed group
 of cardinality $<\chi$) and e.g. G.C.H. (see below).

\noindent
4)  We can improve a little Eklof's results on compactness \cite{Ek82} where
``$A$ free'' is replace by ``${\rm Ext}(A,\bbZ)=0$''. Note that in his proofs
$\Diamond_S$ can be replaced by ``$S$ not small'' e.g.  (see \cite{DvSh:65}),
and instead ``$\Diamond_S$ for stationary $S$'' by the above ``$S$ not
small for all stationary $S$ such that $(\forall \delta \in S) 
\cf(\delta) = \aleph_0$'' suffice but if $\sup (S)= \lambda^+, 
\lambda^{\aleph_0} = \lambda,2^\lambda=\lambda^+$,
this holds. So we can get compactness for $\beth_{\alpha+\omega}$
assuming G.C.H.

\noindent
4A)  Hypothesis I can be rephrased similary to Hypothesis II, as the
existence of a winning strategy (to player II) in appropriate game.
\item[(5)] For the Whithead problem we need only ``any $\lambda$-free
abelian group is $\lambda^+$-free'' for singular $\lambda$. So
suppose $G$ is a $\lambda-$free group with universe $\lambda$
and ${\bf F}=\{(A, B): A/B$ is free\}. There we do not need Hypothesis I, and
can represent the proof somewhat differently.

In the construction we choose \underline{pure} subgroups $A^n_i,
B^n_i$ and choose a free basis $I^n_i$ of $A^n_i$ and demand satisfying
\mn
\begin{enumerate}
\item[$(a)$]  (1) + (2)
\sn
\item[$(b)$]  for $m < n,\, A^m_{i+1} \cap B^n_i$ is generated by a subset of
$I^m_{i+1}$
\sn
\item[$(c)$]  for each $m<n$ and integer $a$,
\[
(\forall x \in B^n_i \cap A^{m+1}_{i+1}) [(\exists y \in A^{m+1}_{i+1})
[a y+x \in A^m_{i+1}] \rightarrow (\exists y \in 
A^{m+1}_{i+1} \cap B^n_i) a y +x \in A^m_{i+1})]
\]
\end{enumerate}
\mn
By (c) we shall get $A^m_{i+1} / 
\bigcup\limits_{m<\omega} B^n_i$ hence it is known (Hill) that

\[
\bigcup\limits_m A^m_{i+1} / \bigcup\limits_m B^m_i = \bigcup\limits_{n}
A^m_i
\]

\mn
is free thus finishing.

A real detailed paper will follow.
\end{remark}
\newpage

\section{On the hypothesis}

\begin{context}
\label{3.0}
${\cU},{\bf F}$ is as in Definition \ref{0.2}.
\end{context}

\begin{notation}
1)  $S_\kappa (A)= \{B \subseteq A :|B|<\kappa\}$.

\noindent
2)  $A/B$ is free mean $(A,B)\in {\bf F}$.

\noindent
3)  $A,B,D$ denote subsets of ${\cU}$.

\noindent
4)  ${\cM} = (\cH(\chi),\in,<^*_\chi)$ where $\chi$ is large enough 
such that ${\cP}({\cU}) \in {\cH}(\chi)$ and $<^*_\chi$ a well
ordering of ${\cM}$.  We say ${\cM}^*$ is a $\kappa$-expansion of 
${\cM}$ if we expand ${\cM}$ by $\le \kappa$ additional relations 
and functions. 

\noindent
5)  ${\cE}^{ub}_\kappa (A)$ is the following filter or
${\cS}_\kappa(A):Y \in \cE^{\ub}_\kappa(A)$ iff 
$Y \supseteq Y_C=Y_C(A)$ for some 
$Y \subseteq A,C \in S_\kappa(A)$ where $Y_C = \{B\in S_\kappa(A):C 
\subseteq B\}$ we call $Y[C,A]$ a generator.
\end{notation}

\begin{definition}
\label{3.3}
1) The pair $A/B$ is ${\cE}$ free $({\cE}$, or ${\cE}(A)$, is a 
filter over a family of subsets of $A$) if

\[
\{C:C \in \cup \cE,C/ B \quad \text{is free}\} \in \cE.
\]

\noindent
2) We can replace ``free'' by any other property.
\end{definition}

\begin{remark}
Obvious monotonicity results hold.
\end{remark}

\begin{definition}
\label{3.4}
For every $\mu \leqq \kappa < \lambda,C \in S_\kappa(A),A$ such that 
$|A|=\lambda$, and $B,$ and filter ${\cE}$ over $S_\kappa(A)$,
we define the rank $R(C,{\cE})$ as an ordinal or $\infty$, so that
\mn
\begin{enumerate}
\item[(a)]  $R(C,{\cE}) \geqq \alpha+1$ iff $C/B$ is free and
$\{D\in S_\kappa(A):C \subseteq D$ and $D/C \cup B$ is free and
$R(D,{\cE}) \geqq \alpha\} \ne \emptyset \mod {\cE}$
\sn
\item[(b)]  $R(C,{\cE}) \geqq \delta (\delta=0$ or $\delta$ limit)
\underline{iff} $C/B$ is free and $\alpha<\delta$ implies 
$R(C,{\cE}) \geqq \alpha$
(more exactly, we should write $R(C,{\cE};A/B)$)
\sn
\item[(c)]  $R(A/ B,{\cE}) = \sup \{R(C,{\cE}):C \in S_\kappa(A)\},
R_\kappa (C) = R_\kappa (C,{\cE}^{ub}_\kappa(A))$
\sn
\item[(d)]  $R^{ub}_\kappa(C) = R(C,{\cE}^{ub}_\kappa)$ and
$R^{ub}_\kappa(A/ B) = R(A/B,{\cE}^{ub}_\kappa)$.
\end{enumerate}
\end{definition}

\begin{lemma}
\label{3.4}
Suppose $\kappa^+<\lambda, \mu \leqq \kappa, A/ B$ is not
${\cE}^{\kappa^+}_{\kappa^+}$-non-free and $S_1 \in 
{\cE}^{\kappa^+}_{\kappa^+}(A)$.

Then $R^{ub}_\kappa=\infty$, [moreover for every
$S_1 \in {\cE}^{\kappa^+}_{\kappa^+}(A)\kappa$-expansion $\cM^*$
of $\cM$ there are $C \in S_2$ and $D \in S_1$ and
$N \prec \cM^*,\{A,B\}\in N,\|N\| = \kappa$ such that 
$D \in N,C=D \cap N$ and $R^\mu_k (C)=\infty$.]
\end{lemma}

\begin{PROOF}{\ref{3.4}}
Let $S_1\supseteq S^M_\kappa (M^*)$ if $C \in S_\kappa(A),0 
\leqq R^{\ub}_\kappa(C)<\infty$, then there is a generator
$S(C)\in {\cE}^{\ub}_\kappa(A),S(C)= S^{\ub}_\kappa(M^*_C)$, such that
for $D \in S(C), D/ C \cup B$ is not free or 
$R^{\ub}_\kappa(D) < R^{\ub}_\kappa(C)$.
If $C\ B$ is not free or $R^{\ub}_\kappa(C) = \infty,$ let 
$\cM^*_C$ be any $\kappa$-expansion of $\cM$,
and let $S_2=S^{\ub}_\kappa(\cM^2)$. Let $\cM^+$ be a
$\kappa$-expansion of $\cM$,expanding 
$\cM^*,\cM^2$ and having the relations $P,P_2$ where

\[
P = \{(C,N):C\in S_\kappa(A),N \prec \cM^*_C,\|N\| <\chi_2\}
\]

\[
P_2 = \{N:N < \cM^2,\|N\| <\chi_2\}.
\]

\mn
As

\[
\{D \in S_{\kappa^+}(A):D/ B \quad \text{ is free }\} \ne \emptyset 
\mod \cE^{\kappa^+}_{\kappa^+} (A)
\]

\mn
and $S_1 \in \cE^{\kappa^+}_{\kappa^+} (A)$ and (by \ref{3.3}
$\cS_{\kappa^+}(A))$; there are $D$, ${\bar N}$ such that:
\mn
\begin{enumerate}
\item[$(1)$]   $S/B$ is free
\sn
\item[$(2)$]  $D \in S_1$
\sn
\item[$(3)$]  $N_i(i<\kappa^+)$ is an
$\underline{M}^+$-sequence and $\|N_i\| \leqq \kappa$, so
\sn
\item[$(4)$]  $D=A \cap \bigcup\limits_{i<\kappa^+} N_i$,
without loss of generality $\|N_i\| = \kappa,\kappa\subseteq N_i$.
\end{enumerate}
\mn
Let $A^*_i = D \cap N_i$, so $A^*_i\in N_{i+1}$ and let
$N=\bigcup\limits_{i<\kappa^+} N_i$. 
Clearly $\langle N_i:i<\kappa^+\rangle$ is also an $\cM^2$-sequence
hence for each $\delta<\kappa^+, \langle N_i:i<\delta\rangle$ is an
$\cM^2$-sequence, hence, if $\kappa$ divides 
$\delta,\cf(\delta)=\mu$, then $A^*_\delta \in S_2$. 
If $C \in N_i,C \in S_\kappa(A)$, then for every $j>i,j<\kappa^+$ there
is a model $N^j_i \prec \cM^*_C,\|N^j_i\| = \kappa,
|N^j_i|$ and $N^i_j\in N_{j+1},$ hence $N^i_j \subseteq N_{j+1}$.

Hence, for any limit ordinal $\delta, i<\delta<\kappa^+$ implies
$N_\delta \prec M^*_C$.  Clearly $\langle N_j:i<j<\kappa^+, j$ limit
$\rangle$ is an $\underline{M}^+$-sequence, hence it is an 
$\underline{M}^*_C$ sequence, hence, is $i<\delta<\kappa^+,\delta$ is 
limit, $\kappa^2$ divides $\delta,\cf(\delta) = \mu$, then
$A^*_\delta\in S(C)$. As $S/ B$ is free, by \cite{Sh:52},1.2(7) there is a
closed unbounded subset of $\kappa^+,W$, such that for 
$i,j \in W,i<j,A^*_j/A^*_i\cup B$ is free and $A^*_i/ B$ is free. 
We can assume that such $i \in W$ is divisible by $\kappa^2$. 
Hence, if $i,j \in W,i<j,\cf(j)=\mu, R^{\ub}_\kappa(A^*_i)<\infty$, then
$R^\mu_\kappa (A^*_j)< R^\mu_\kappa (A^*_i)<\infty$
(by the definition of $S(C))$.
So, if for some $i\in W, R^\mu_\kappa(A^*_i)<\infty,\cf(i_n) = \mu,i_n\in
W,i<i_n<i_{n+1}$ then $R^\mu_\kappa(A^*_{i_n})$ is an 
infinite decreasing sequence of ordinals, a contradiction. 
Hence, $i\in W$ implies $R^\mu_\kappa(A^*_i)=\infty$. Let $D=
\bigcup\limits_{i<\kappa^+} A^*_i$, and choose $N \prec
\underline{M}^*,D \in N,N \cap \bigcup\limits_{i<\kappa^+}
A^*_i= A^*_\delta,\delta \in W, \cf(\delta) =\mu$, and
$C=A^*_\delta$. So we are finished.
\end{PROOF}

\begin{lemma}
\label{3.5}
1) If $\mu \leqq \kappa < \lambda, C\in S_\kappa(A), R^\mu_\kappa(C)
=\infty, S\in {\cE}^\mu_\kappa (A)$, then for some $D \in S, 
C \subseteq D, R^\mu_\kappa(D)=\infty$ and $D/C\cup B$ is free.

\noindent
2)  The same holds for any filter over $S_\kappa(A)$.
\end{lemma}

\begin{PROOF}{\ref{3.5}}
1)  As $S_\kappa(A)$ is a set, for some ordinal $\alpha_0 < 
|S_\kappa(A)|^+$, for no $C \in S_\kappa(A)$ is 
$R^\mu_\kappa(C)=\alpha_0$. We can easily prove that
$R^\mu_\kappa(C) \geqq \alpha_0$ iff $R^\mu_\kappa(C)=\infty$.
Using the definition we get our assertion.

\noindent
2)  The same proof.
\end{PROOF}
\newpage

\bibliographystyle{amsalpha}
\bibliography{lista,listb,listx,listf,liste,listz}

\end{document}